\documentclass{tran-l}

\title[Darcy's law for a compressible thermofluid]
{Darcy's law for a compressible thermofluid}

\author{Anvarbek Meirmanov}

\date{}
\newtheorem{theorem}{Theorem}[section]
\newtheorem{lemma}{Lemma}[section]

\newtheorem{corollary}{Corollary}[section]
\newtheorem{definition}{Definition}[section]
\newtheorem{assumption}{Assumption}[section]

\renewcommand{\div}{\mbox{div}}
\newcommand{\x}{{\mathbf x}}
\newcommand{\y}{{\mathbf y}}

\newcommand{\w}{{\mathbf w}}

\newcommand{\uu}{{\mathbf u}}
\newcommand{\vv}{{\mathbf v}}
\newcommand{\e}{{\mathbf e}}
\newcommand{\n}{{\mathbf n}}
\renewcommand{\mathbf}[1]{\mbox{\boldmath$#1$}}
\newcommand{\V}{\mathbf V}
\newcommand{\W}{\mathbf W}

\newcommand{\D}{\mathbb D}

\newcommand{\RR}{\mathbb R}
\newcommand{\PP}{\mathbb P}

\newcommand{\I}{\mathbb I}

\begin{document}

\maketitle \noindent \textbf{Abstract.} A linear system of
differential equations describing a joint motion of a thermoelastic
porous body with a sufficiently large Lam\'{e}'s  constants
(absolutelty rigid body) and a thermofluid, occupying porous space,
is considered.
 The rigorous justification, under various
conditions imposed on physical parameters, is fulfilled for
homogenization procedures as the dimensionless size of the pores
tends to zero, while the porous body is geometrically periodic. As
the results we derive Darcy's system of filtration for thermofluid,
depending on ratios between physical parameters. The proofs are
based on Nguetseng's two-scale
convergence method of homogenization in periodic structures.\\

\noindent \textbf{Key words:}  Anisothermic Stokes and Lam\'{e}'s
 equations, two-scale convergence, homogenization of
periodic structures.\\
\noindent \textbf{MOS subject classification:} 35M99;76Q05\\

\address{State University, ul. Pobedi, 85, Belgorog,
 Russia, email:meirmanov@bsu.edu.ru}\\
\normalsize

\addtocounter{section}{0} \setcounter{equation}{0}

\begin{center} \textbf{Introduction}
\end{center}
In the present publication we consider a problem of a joint motion
of thermoelastic deformable solid (thermoelastic skeleton),
perforated by a system of channels (pores) and incompressible
thermofluid occupying a porous space. In dimensionless variables
(without primes)

$$ \x'=L \x,\quad t'=\tau t,\quad \w'=L \w,
 \quad \theta'=\vartheta_*\frac{L }{\tau v_{*}} \theta,\quad \rho'_s=
\rho_0 \rho_s,\quad \rho'_f =\rho_0 \rho_f,$$

the differential equations of the problem in a domain $\Omega \in
\RR^{3}$ for the dimensionless displacement vector $\w$ of the
continuum medium and the dimensionless temperature $\theta$, have a
form:
\begin{eqnarray} \label{0.1}
& \displaystyle \alpha_\tau \bar{\rho} \frac{\partial^2
\w}{\partial t^2}=\div_x \PP,\\
 \label{0.2}
& \displaystyle \alpha_\tau \bar{c}_p \frac{\partial
\theta}{\partial t} = \div_x ( \bar{\alpha} _{\kappa} \nabla_x
\theta) -\bar{\alpha}_\theta
\frac{\partial}{\partial t} \div_x \w,\\
 \label{0.3} & \displaystyle \PP =
\bar{\chi}\alpha_\mu \D\Bigl(\x,\frac{\partial \w}{\partial t}\Bigr)
+(1-\bar{\chi})\alpha_\lambda \D(x,\w)-(q+\alpha_{\theta s}(1-\bar{\chi})\theta)\I ,\\
 \label{0.4}
& \displaystyle q=p+\frac{\alpha_\nu}{\alpha_p}\frac{\partial
p}{\partial t}+\alpha_{\theta f}\bar{\chi} \theta,\\
 \label{0.5}
& \displaystyle p+\bar{\chi} \alpha_p \div_x \w=0.
\end{eqnarray}

 The problem is endowed   with initial and boundary conditions
\begin{equation} \label{0.6}
\w|_{t=0}=0,\quad \frac{\partial \w}{\partial t}|_{t=0}=0,\quad
\theta|_{t=0} =0,\quad \x\in \Omega
\end{equation}
\begin{equation} \label{0.7}
\w=\w_0,\quad \theta=\theta_0,\quad \x \in S=\partial \Omega, \quad
t\geq 0.
\end{equation}
Here and further we use  notations
 $$ \D(x,\uu)=(1/2)\left(\nabla_x \uu +(\nabla_x \uu)^T\right),$$
$$\bar{\rho}=\bar{\chi}\rho_f +(1-\bar{\chi})\rho_s, \quad
\bar{c}_p=\bar{\chi} c_{pf} +(1-\bar{\chi})c_{ps},$$
$$  \bar{\alpha _{\kappa}} =\bar{\chi} \alpha _{\kappa f}
+(1-\bar{\chi})\alpha _{\kappa s},\quad \bar{\alpha}_\theta
=\bar{\chi} \alpha_{\theta f} +(1-\bar{\chi})\alpha_{\theta s}.$$ In
this model the characteristic function of the porous space
$\bar{\chi}(\x)$ is a known function.

 For more details about  Eqs.\eqref{0.1}-- \eqref{0.5} and description of
 dimensionless constants (all these constants are positive ) see \cite{MS}.

We accept the following constraints
\begin{assumption} \label{assumption1}
domain  $\Omega =(0,1)^3$ is a periodic repetition of an elementary
cell  $Y^\varepsilon =\varepsilon Y$, where $Y=(0,1)^3$ and quantity
$1/\varepsilon$ is integer, so that $\Omega$ always contains an
integer number of elementary cells $Y_i^\varepsilon$. Let $Y_s$  be
a "solid part" of $Y$, and the "liquid part"  $Y_f$ -- is its open
complement. We denote as $\gamma = \partial Y_f \cap \partial Y_s$
and $\gamma $ is $C^{1}$-surface.
 A porous space  $\Omega ^{\varepsilon}_{f}$  is the periodic repetition of
the elementary cell $\varepsilon Y_f$, and solid skeleton  $\Omega
^{\varepsilon}_{s}$  is the periodic repetition of the elementary
cell $\varepsilon Y_s$. A boundary  $\Gamma^\varepsilon =\partial
\Omega_s^\varepsilon \cap \partial \Omega_f^\varepsilon$  is the
periodic repetition in  $\Omega$ of the boundary $\varepsilon
\gamma$. The "solid skeleton"  $\Omega _{s}$ is a connected domain.
\end{assumption}
In these assumptions
\begin{equation*}
 \bar{\chi}(\x)=\chi^{\varepsilon}(\x)=\chi
 \left(\x / \varepsilon\right),
\end{equation*}
$$\bar{c}_{p}=c_{p}^{\varepsilon}(\x)=\chi^{\varepsilon}(\x)c _{pf}+
(1-\chi^{\varepsilon}(\x))c_{ps},$$
$$\bar{\rho}=\rho^{\varepsilon}(\x)=\chi^{\varepsilon}(\x)\rho _{f}+
(1-\chi^{\varepsilon}(\x))\rho_{s},$$
$$  \bar{\alpha} _{\kappa} =\alpha^{\varepsilon} _{\kappa}(\x)=
\chi ^{\varepsilon}(\x)\alpha _{\kappa f} +(1-\chi
^{\varepsilon}(\x))\alpha _{\kappa s}, $$ $$\bar{\alpha}_\theta
=\alpha ^{\varepsilon}_\theta(\x)=\chi ^{\varepsilon}(\x)
\alpha_{\theta f} +(1-\chi ^{\varepsilon}(\x))\alpha_{\theta s},$$

Let $\varepsilon$ be a characteristic size of pores $l$ divided by
the characteristic size $L$ of the entire porous body:
$$\varepsilon =\frac{l}{L}.$$

Suppose that all dimensionless parameters depend on the small
parameter $\varepsilon$ and there exist limits (finite or infinite)
$$\lim_{\varepsilon\searrow 0} \alpha_\mu(\varepsilon) =\mu_0, \quad
\lim_{\varepsilon\searrow 0} \alpha_\lambda(\varepsilon) =\lambda_0,
\quad \lim_{\varepsilon\searrow 0}
\alpha_\tau(\varepsilon)=\tau_{0},$$

$$\lim_{\varepsilon\searrow 0} \alpha_\eta(\varepsilon) =\eta_0,
\quad \lim_{\varepsilon\searrow 0} \alpha_p(\varepsilon) =p_{*},
\quad \lim_{\varepsilon\searrow 0}\alpha_\nu(\varepsilon) =\nu_0,$$

$$\lim_{\varepsilon\searrow 0} \alpha _{\kappa
f}=\kappa_{0f},\quad \lim_{\varepsilon\searrow 0} \alpha_{\kappa
s}(\varepsilon) =\kappa _{0s},\quad \lim_{\varepsilon\searrow 0}
\alpha _{ \theta f}(\varepsilon) =\beta_{0f},$$

$$ \lim_{\varepsilon\searrow 0} \alpha _{ \theta s}(\varepsilon)
=\beta_{0s},\quad \lim_{\varepsilon\searrow 0}
\frac{\alpha_\mu}{\varepsilon^{2}} =\mu_1.$$
 We restrict our consideration by the case when when
 $$\mu_0=0; \quad \lambda_0=\eta_0= \infty ; \quad   \beta_{0s},
 \,\beta_{0f},\,\nu_0\,\tau_{0} <\infty,$$
$$ 0<\mu_1+\tau_{0}, \, p_{*},\,\kappa_{0f},\,\kappa_{0s} <\infty.$$

The condition $\lambda_0=\eta_0= \infty $  means that the solid
skeleton is an absolutely rigid body.

 Using  Nguetseng's two-scale convergence method \cite{NGU,LNW} we show that the
 limiting regime as $\varepsilon$ tends to zero is described by two types of Darcy's like
 system of equations of filtration for the
velocity of the liquid component, coupled with corresponding heat
equation.

Different isothermic  models have been considered in
 \cite{B-K}--\cite{AM}.
 The first research with the aim of finding limiting regimes in the
case when the skeleton was assumed to be an absolutely rigid
isothermal body was carried out by E. Sanchez-Palencia and L.
Tartar. E. Sanchez-Palencia \cite[Sec. 7.2]{S-P} formally obtained
Darcy's law of filtration  using the method of  two-scale asymptotic
expansions, and L. Tartar \cite[Appendix]{S-P} mathematically
rigorously justified the homogenization procedure. Using the same
method of  two-scale expansions J. Keller and R. Burridge \cite{B-K}
derived formally the system of Biot's equations
 in the case when
the parameter $\alpha _{\mu}$  was of order $\varepsilon^2$, and the
rest of the coefficients were fixed independent of $\varepsilon$.
 Under the same assumptions as in the article \cite{B-K}, the rigorous
justification of Biot's model has been given by G. Nguetseng
\cite{GNG} and later by A. Mikeli\'{c}, R. P. Gilbert, Th. Clopeaut,
and J. L. Ferrin in \cite{G-M2,G-M3,G-M1}.  The most general case
has been studied in \cite{AM}.

\addtocounter{section}{1}
\setcounter{theorem}{0} \setcounter{lemma}{0}
\setcounter{proposition}{0} \setcounter{corollary}{0}
\setcounter{definition}{0} \setcounter{assumption}{0}

\begin{center} \textbf{\S1. Formulation of the main results.}
\end{center}

As usual, equations \eqref{0.1}-\eqref{0.2} are understood in the
sense of distributions. They involve the equations \eqref{0.1}--
\eqref{0.5} in a usual sense in the domains $\Omega_f^{\varepsilon}$
and $\Omega_s^{\varepsilon}$ and the boundary conditions
\begin{eqnarray} \label{1.1}
& [\vartheta]=0, \quad [\w]=0,\quad \x_0\in \Gamma ^{\varepsilon},\; t\geq 0,\\
\label{1.2} & [\PP]=0,\quad [\alpha ^{\varepsilon} _{\kappa}
\nabla_x \theta ]=0, \quad \x_0\in \Gamma ^{\varepsilon},\; t\geq 0
\end{eqnarray}
on the boundary  $\Gamma^\varepsilon $, where
\begin{eqnarray}
\nonumber & [\varphi](\x_0)=\varphi_{(s)}(\x_0)
-\varphi_{(f)}(\x_0),\\
 \nonumber \displaystyle
& \varphi_{(s)}(\x_0) =\lim\limits_{\tiny \begin{array}{l}\x\to \x_0\\
\x\in \Omega_s^{\varepsilon}\end{array}} \varphi(\x),\quad
\varphi_{(f)}(\x_0) =\lim\limits_{\tiny \begin{array}{l}\x\to \x_0\\
\x\in \Omega_f^{\varepsilon}\end{array}} \varphi(\x).
\end{eqnarray}

 There are various equivalent in the sense of
distributions forms of representation of equations
\eqref{0.1}--\eqref{0.2} and boundary conditions
\eqref{1.1}--\eqref{1.2}. In what follows, it is convenient to write
them in the form of the integral equalities.

\begin{definition} \label{definition1}
Functions
$(\w^{\varepsilon},\theta^{\varepsilon},p^{\varepsilon},q^{\varepsilon}$
are called a generalized solution of the problem \eqref{0.1}--
\eqref{0.7}, if they satisfy the regularity conditions in the domain
$ \Omega_{T}=\Omega\times (0,T)$
\begin{equation} \label{1.3}
\w^{\varepsilon},\, \D(x,\w^{\varepsilon}),\,
\div_x\w^{\varepsilon},\, q^{\varepsilon},\,p^{\varepsilon},\,
\frac{\partial p^{\varepsilon}}{\partial
 t},\,\theta^{\varepsilon}, \nabla_x \theta
^{\varepsilon} \in L^2(\Omega_{T})
\end{equation}
in the domain $\Omega_{T}=\Omega\times (0,T)$, boundary conditions
\eqref{0.7} with functions
$$\w^{\varepsilon}_{0},\, \theta^{\varepsilon}_{0} \in  L^2
((0,T);W^1_2(\Omega )),$$
 equations
\begin{eqnarray} \label{1.4}
&\displaystyle q^{\varepsilon}=p^{\varepsilon}+
\frac{\alpha_\nu}{\alpha_p}\frac{\partial p^{\varepsilon}}{\partial
t}+\alpha_{\theta f}\chi^{\varepsilon}\theta ^{\varepsilon},\\
\label{1.5}& \displaystyle p^{\varepsilon}+ \chi^{\varepsilon}
\alpha_p \div_x \w^{\varepsilon}=0
\end{eqnarray}
 a.e. in  $\Omega_{T}$, integral identity
\begin{eqnarray}\label{1.6}
\left. \begin{array}{lll} \displaystyle\int_{\Omega_{T}}
\Bigl(\alpha_\tau \rho ^{\varepsilon} \w^{\varepsilon}\cdot
\frac{\partial ^{2}{\mathbf \varphi}}{\partial t^{2}} - \chi
^{\varepsilon}\alpha_\mu \D(\x, \w^{\varepsilon}):
\D(x,\frac{\partial {\mathbf \varphi}}{\partial t})+\\[1ex]
\{(1-\chi^{\varepsilon})\alpha_\lambda\D(x,\w^{\varepsilon})-
(q^{\varepsilon}+\alpha_{\theta s}(1-\chi^{\varepsilon})\theta
^{\varepsilon})\I\} : \D(x,{\mathbf \varphi})\Bigr) d\x dt=0
\end{array} \right\}
\end{eqnarray}
for all smooth  vector-functions  ${\mathbf \varphi}={\mathbf
\varphi}(\x,t)$ such that  ${\mathbf \varphi}|_{\partial \Omega}
={\mathbf \varphi}|_{t=T}=\partial {\mathbf \varphi} / \partial
t|_{t=T}=0$ and integral identity
\begin{eqnarray}\label{1.7}
 \int_{\Omega_{T}} \Bigl((\alpha_\tau c^{\varepsilon}_p
\theta ^{\varepsilon}+\alpha^{\varepsilon}_\theta \div_x \w
^{\varepsilon}) \frac{\partial \xi}{\partial t} - \alpha _{\kappa
}^{\varepsilon} \nabla_x \theta ^{\varepsilon}\cdot \nabla_x \xi
+\Psi \xi \Bigr) d\x dt=0
\end{eqnarray}
 for all smooth functions $\xi= \xi(\x,t)$ such that $\xi|_{\partial \Omega}
= \xi|_{t=T}=0$.
\end{definition}

In \eqref{1.6}  by $A:B$ we denote the convolution (or,
equivalently, the inner tensor product) of two second-rank tensors
along the both indexes, i.e., $A:B=\mbox{tr\,} (B^*\circ
A)=\sum_{i,j=1}^3 A_{ij} B_{ji}$.

We suppose the next assumption to be held:
\begin{assumption} \label{assumption2}
 Sequences  $\{(\sqrt{\alpha_\mu}+\sqrt{\alpha_\nu})\chi ^{\varepsilon}|\nabla
\partial \w^{\varepsilon}_0 /\partial t|\}$, $\{\sqrt{\alpha_\lambda}(1-\chi
^{\varepsilon})\nabla \partial \w^{\varepsilon}_0 /
\partial t\}$, $\{(\sqrt{\alpha _{p}}+\sqrt{\alpha _\theta
^{\varepsilon}}) \div_x\partial \w^{\varepsilon}_0 /
\partial t\}$, $\{\sqrt{\alpha_\tau}\partial^{2} \w^{\varepsilon}_0 /
\partial t^{2}\}$, $\{\sqrt{\alpha_\tau}\partial
\theta^{\varepsilon}_0/\partial t\}$, $\{\nabla
 \theta^{\varepsilon}_0\}$  are uniformly in $\varepsilon$  bounded in
$L^2(\Omega).$
\end{assumption}

  In what follows all parameters may take all permitted values.  If, for
example, $\tau_{0}=0$ , then all terms in final equations
 containing this  parameter  disappear.

 The following Theorems \ref{theorem1}--\ref{theorem2} are the main results of the paper.

\begin{theorem} \label{theorem1}
For all $\varepsilon >0$  on the arbitrary time interval  $[0,T]$
there exists a unique generalized solution of the problem
\eqref{0.1}-- \eqref{0.7}  and
\begin{equation} \label{1.8}
 \displaystyle \max\limits_{0\leq t\leq
T}\| |\w^{\varepsilon}(t)|+\sqrt{\alpha_\mu} \chi^\varepsilon
|\nabla_x \w^{\varepsilon}(t)|+(1-\chi^\varepsilon)
 \sqrt{\alpha_\lambda}|\nabla_x \w^{\varepsilon}(t)| \|_{2,\Omega}
   \leq C_{0} ,
\end{equation}
\begin{equation} \label{1.9}
 \displaystyle\| \theta^{\varepsilon} \|_{2,\Omega_{T}}+
 \|  \nabla_x \theta^{\varepsilon}\|_{2,\Omega_{T}} \leq C_{0} ,
\end{equation}
\begin{equation}\label{1.10}
 \|q^{\varepsilon}\|_{2,\Omega_{T}} +
\|p^{\varepsilon}\|_{2,\Omega_{T}} + \frac{\alpha _{\nu}}{\alpha
_{p}}\|\frac{\partial p^{\varepsilon}}{\partial
t}\|_{2,\Omega_{T}}\leq C_{0}
\end{equation}
where $C_{0}$ does not depend on the small parameter $\varepsilon $.
\end{theorem}

\begin{theorem} \label{theorem2}
Sequences
 $\{\chi ^{\varepsilon}\w^\varepsilon\}$,
 $\{p^{\varepsilon}\}$  and $\{q^{\varepsilon}\}$ converge weakly
  in  $L^{2}(\Omega_{T})$ to  $\w ^{f}$,  $p$,  and  $q$
  respectively.  The sequence $\{\theta^{\varepsilon}\}$ converge weakly
  in $L^{2}((0,T);W^1_2(\Omega))$  to function  $\theta$.  Functions $\w^{\varepsilon}$
   admit an extension  $\uu^{\varepsilon}$  from $\Omega_{s}^{\varepsilon}\times (0,T)$
 into $\Omega_{T}$ such that the sequence $\{\uu^{\varepsilon}\}$  converge strongly in
 $L^{2}((0,T);W^1_2(\Omega))$ to zero. Weak and strong limits  $\vv=\partial \w ^{f}/\partial t$,
  $q$, $p$ and  $\theta$  satisfy in the domain  $\Omega_{T}$ the state equation
\begin{equation}\label{1.11}
  q=p +\nu_0 p_{*}^{-1} \partial p/\partial t
+m\beta_{0f}\theta ,
\end{equation}
the continuity equation
\begin{equation}\label{1.12}
 (1/p_{*})p+\div_x \w =0 ,
\end{equation}
the heat equation
\begin{equation}\label{1.13}
\tau_{0}\hat{c_{p}}\frac{\partial \theta}{\partial t}
-\frac{\beta_{0f}}{p_{*}}\frac{\partial p}{\partial t}
 =\div_x ( B^{\theta}\cdot \nabla \vartheta )
\end{equation}
and  Darcy's law in the form
\begin{equation}\label{1.14}
\vv=-\int_{0}^{t} B_{1}(\mu_1,t-\tau)\cdot \nabla_x q(\x,\tau )d\tau
\end{equation}
 in the case of $\tau_{0}>0$ and
$\mu_{1}>0$;  Darcy's law in the form
\begin{equation}\label{1.15}
\vv=-B_{2}(\mu_1)\cdot\nabla_x q
\end{equation}
in the case of $\tau_{0}=0$  and, finally, Darcy's law in the form
\begin{equation}\label{1.16}
\vv=-\frac{1}{\tau_{0}\rho_{f}}(m\I-B_{3})\cdot\int_{0}^{t}\nabla_x
q(\x,\tau )d\tau
\end{equation}
 in the case of $\mu_{1}=0$. The
problem is supplemented by  the boundary condition
\begin{equation}\label{1.17}
(\vv-\vv_{0})\cdot \n(\x)=0, \quad \x \in S, \quad t>0
\end{equation}
for the velocity of the liquid component and boundary and initial
conditions
\begin{equation}\label{1.18}
 \tau _{0}\theta(\x,0)=0,\quad \theta(\x,0)=\theta _{0}(\x)
\end{equation}
for the temperature.

In Eq.\eqref{1.13}  symmetric strictly positively defined matrix
$B^{\theta}$ is given  below by formula \eqref{4.19} and in Eqs.
\eqref{1.14}--\eqref{1.16} $\n(\x)$ is the unit normal vector to $S$
at a point $\x\in S$, and symmetric strictly positively defined
matrices $B_{1}(\mu_1,t)$, $B_{2}(\mu_1)$, and $B_{3}$ are defined
in \cite{AM}.
\end{theorem}

 \addtocounter{section}{1}
\setcounter{theorem}{0} \setcounter{lemma}{0}
\setcounter{proposition}{0} \setcounter{corollary}{0}
\setcounter{definition}{0} \setcounter{assumption}{0}

\begin{center} \textbf{\S2. Preliminaries}
\end{center}

\textbf{2.1. Two-scale convergence.} Justification of Theorems
\ref{theorem1}--\ref{theorem2} relies on systematic use of the
method of two-scale convergence, which had been proposed by G.
Nguetseng \cite{NGU} and has been applied recently to a wide range
of homogenization problems (see, for example, the survey
\cite{LNW}).

\begin{definition} \label{TS}
A sequence $\{\varphi^\varepsilon\}\subset L^2(\Omega_{T})$ is said
to be \textit{two-scale convergent} to a limit $\varphi\in
L^2(\Omega_{T}\times Y)$ if and only if for any 1-periodic in $\y$
function $\sigma=\sigma(\x,t,\y)$ the limiting relation
\begin{equation}\label{(2.1)}
\lim_{\varepsilon\searrow 0} \int_{\Omega_{T}}
\varphi^\varepsilon(\x,t) \sigma\left(\x,t,\x /
\varepsilon\right)d\x dt = \int _{\Omega_{T}}\int_Y
\varphi(\x,t,\y)\sigma(\x,t,\y)d\y d\x dt
\end{equation}
holds.
\end{definition}

Existence and main properties of weakly convergent sequences are
established by the following fundamental theorem \cite{NGU,LNW}:
\begin{theorem} \label{theorem3}(\textbf{Nguetseng's theorem})

\textbf{1.} Any bounded in $L^2(Q)$ sequence contains a subsequence,
two-scale convergent to some limit
$\varphi\in L^2(\Omega_{T}\times Y)$.\\[1ex]
\textbf{2.} Let sequences $\{\varphi^\varepsilon\}$ and
$\{\varepsilon \nabla_x \varphi^\varepsilon\}$ be uniformly bounded
in $L^2(\Omega_{T})$. Then there exist a 1-periodic in $\y$ function
$\varphi=\varphi(\x,t,\y)$ and a subsequence
$\{\varphi^\varepsilon\}$ such that $\varphi,\nabla_y \varphi\in
L^2(\Omega_{T}\times Y)$, and $\varphi^\varepsilon$ and $\varepsilon
\nabla_x \varphi^\varepsilon$ two-scale converge to $\varphi$ and
$\nabla_y \varphi$,
respectively.\\[1ex]
\textbf{3.} Let sequences $\{\varphi^\varepsilon\}$ and $\{\nabla_x
\varphi^\varepsilon\}$ be bounded in $L^2(Q)$. Then there exist
functions $\varphi\in L^2(\Omega_{T})$ and $\psi \in
L^2(\Omega_{T}\times Y)$ and a subsequence from
$\{\varphi^\varepsilon\}$ such that $\psi$ is 1-periodic in $\y$,
$\nabla_y \psi\in L^2(\Omega_{T}\times Y)$, and
$\varphi^\varepsilon$ and $\nabla_x \varphi^\varepsilon$ two-scale
converge to $\varphi$ and $\nabla_x \varphi(\x,t)+\nabla_y
\psi(\x,t,\y)$, respectively.
\end{theorem}

\begin{corollary} \label{corollary2.1}
Let $\sigma\in L^2(Y)$ and
$\sigma^\varepsilon(\x):=\sigma(\x/\varepsilon)$. Assume that a
sequence $\{\varphi^\varepsilon\}\subset L^2(\Omega_{T})$ two-scale
converges to $\varphi \in L^2(\Omega_{T}\times Y)$. Then the
sequence $\sigma^\varepsilon \varphi^\varepsilon$ two-scale
converges to $\sigma \varphi$.
\end{corollary}

\textbf{2.2. An extension lemma.} The typical difficulty in
homogenization problems while passing to a limit in Model
${(\mathbf{N}\mathbf B})^\varepsilon$  as $\varepsilon \searrow 0$
arises because of the fact that the bounds on the gradient of
displacement $\nabla_x \w^\varepsilon$ may be distinct in liquid and
rigid phases. The classical approach in overcoming this difficulty
consists of constructing of extension to the whole $\Omega$ of the
displacement field defined merely on $\Omega_s$. The following lemma
is valid due to the well-known results from \cite{ACE,JKO}. We
formulate it in appropriate for us form:

\begin{lemma} \label{Lemma2.1}
Suppose that assumptions of Sec. 1.2 on geometry of periodic
structure hold,  $ \psi^\varepsilon\in W^1_2(\Omega^\varepsilon_s)$
and   $\psi^\varepsilon =0$ on $S_{s}^{\varepsilon}=\partial \Omega
^\varepsilon_s \cap
\partial \Omega$ in the trace sense.  Then there exists a function
$ \sigma^\varepsilon \in
 W^1_2(\Omega)$ such that its restriction on the sub-domain
$\Omega^\varepsilon_s$ coincide with $\psi^\varepsilon$, i.e.,
\begin{equation} \label{2.2}
(1-\chi^\varepsilon(\x))( \sigma^\varepsilon(\x) - \psi^\varepsilon
(\x))=0,\quad \x\in\Omega,
\end{equation}
and, moreover, the estimate
\begin{equation} \label{2.3}
\|\sigma^\varepsilon\|_{2,\Omega}\leq C\|
\psi^\varepsilon\|_{2,\Omega ^{\varepsilon}_{s}}  , \quad \|\nabla_x
\sigma^\varepsilon\|_{2,\Omega} \leq  C \|\nabla_x
 \psi^\varepsilon\|_{2,\Omega ^{\varepsilon}_{s}}
\end{equation}
 hold true, where the constant $C$
depends only on geometry $Y$ and does not depend on $\varepsilon$.
\end{lemma}

\textbf{2.3. Friedrichs--Poincar\'{e}'s inequality in periodic
structure.} The following lemma was proved by L. Tartar in
\cite[Appendix]{S-P}. It specifies Friedrichs--Poincar\'{e}'s
inequality for $\varepsilon$-periodic structure.
\begin{lemma} \label{F-P}
Suppose that assumptions on the geometry of $\Omega^\varepsilon_f$
hold true. Then for any function $\varphi\in
\stackrel{\!\!\circ}{W^1_2}(\Omega^\varepsilon_f)$ the inequality
\begin{equation} \label{(F-P)}
\int_{\Omega^\varepsilon_f} |\varphi|^2 d\x \leq C \varepsilon^2
\int_{\Omega^\varepsilon_f} |\nabla_x \varphi|^2 d\x
\end{equation}
holds true with some constant $C$, independent of $\varepsilon$.
\end{lemma}

\textbf{2.4. Some notation.} Further we denote

 1) $ \langle\Phi \rangle_{Y} =\int_Y \Phi  dy, \quad
 \langle\Phi \rangle_{Y_{f}} =\int_{Y_{f}} \Phi  dy,\quad
 \langle\Phi \rangle_{Y_{s}} =\int_{Y_{s}} \Phi  dy.$

2) If $\textbf{a}$ and $\textbf{b}$ are two vectors then the matrix
$\textbf{a}\otimes \textbf{b}$ is defined by the formula
$$(\textbf{a}\otimes \textbf{b})\cdot
\textbf{c}=\textbf{a}(\textbf{b}\cdot \textbf{c})$$ for any vector
$\textbf{c}$.

\addtocounter{section}{1}
\setcounter{theorem}{0} \setcounter{lemma}{0}
\setcounter{proposition}{0} \setcounter{corollary}{0}
\setcounter{definition}{0} \setcounter{assumption}{0}

\begin{center} \textbf{\S3. Proof of Theorem  \ref{theorem1}}
\end{center}

  Under restriction $\tau_{0}>0$ estimates
 \eqref{1.8}-\eqref{1.9} follow from estimates
\begin{equation*}
\max\limits_{0<t<T}(\sqrt{\alpha_\lambda}\|(1-\chi
^{\varepsilon})\nabla_x \w^{\varepsilon}(t) \|_{2,\Omega} +
\sqrt{\alpha_\tau}\| \frac{\partial \w^\varepsilon}{\partial
t}(t)\|_{2,\Omega}+
\end{equation*}
\begin{equation*}
\sqrt{\alpha _{p}} \|\chi ^{\varepsilon}\div_x
\w^{\varepsilon}(t)\|_{2,\Omega}+
\sqrt{\alpha_\tau}\|\theta^{\varepsilon}(t)\|_{2,\Omega}) +
\|\nabla_x \theta^{\varepsilon}\|_{2,\Omega _{T}}+
\end{equation*}
\begin{equation} \label{3.1}
\sqrt{\alpha_\mu}\|\chi ^{\varepsilon}
 \nabla_x \frac{\partial\w^\varepsilon}{\partial t} \|_{2,\Omega_T}+
 \sqrt{\alpha _{\nu}}\| \chi ^{\varepsilon} \div_x
\frac{\partial\w^\varepsilon}{\partial t}\|_{2,\Omega _{T}} \leq
C_{0},
\end{equation}
where  $C_{0}$ is independent of  $\varepsilon$. These estimates we
obtain if we multiply equation for  $\w^{\varepsilon}$ by $\partial
\w^{\varepsilon}/\partial t$, equation for  $\theta^{\varepsilon}$
  multiply  by  $\theta^{\varepsilon}$, integrate by parts and sum the result.
 The same estimates \eqref{3.1} guaranties the existence and uniqueness of the generalized solution
for the problem \eqref{0.1}--\eqref{0.7}.

Estimates \eqref{1.10}  for pressures follows from
Eqs.\eqref{1.4}--\eqref{1.5} and estimates \eqref{3.1}.

 Estimation of $\w^\varepsilon$ and $\theta^\varepsilon$ in the
case $\tau_0=0$ is not simple, and we outline it in more detail.

First of all we use estimates \eqref{3.1} in the form
\begin{equation}\label{3.2}
\sqrt{\alpha_\mu}\|\chi ^{\varepsilon}
 \nabla_x\w^{\varepsilon}(t)\|_{2,\Omega}+
 \sqrt{\alpha_\lambda}\|(1- \chi ^{\varepsilon})
\nabla_x\w^{\varepsilon}(t) \|_{2,\Omega } \leq C_{0}.
\end{equation}
Next, on the strength of Lemma \ref{Lemma2.1}, we construct
extensions $\uu^\varepsilon $ and $\uu_{0}^\varepsilon $ of the
functions  $\w^\varepsilon $ and $\w_{0}^\varepsilon$  from
$\Omega_s^\varepsilon$ into $\Omega_f^\varepsilon$, such that
$\uu^\varepsilon=\w^\varepsilon$ and $\uu_{0}^\varepsilon
=\w_{0}^\varepsilon$ in $\Omega_s^\varepsilon$,
$\uu_{0}^\varepsilon$ uniformly with respect to $\varepsilon$
bounded in $W_2^1(\Omega)$ and
$$\| \uu^\varepsilon-\uu_{0}^\varepsilon \|_{2,\Omega} \leq C
\|\nabla_x(\uu^\varepsilon-\uu_{0}^\varepsilon )\|_{2,\Omega} \leq
\frac{C}{\sqrt{\alpha_\lambda}}
 \|(1-\chi^\varepsilon)\sqrt{\alpha_\lambda}\nabla_x
 (\w^\varepsilon-\w_{0}^\varepsilon)\|_{2,\Omega }.$$

After that we estimate $\|\w^\varepsilon\|_{2,\Omega}$ with the help
of  Friedrichs--Poincar\'{e}'s inequality in periodic structure
(lemma \ref{F-P}) for the difference $(\uu^\varepsilon
-\w^\varepsilon)$  and estimates \eqref{3.2}:

$$\|\w^\varepsilon\|_{2,\Omega} \leq
\|\uu^\varepsilon\|_{2,\Omega} + \|\uu^\varepsilon
-\w^\varepsilon\|_{2,\Omega} \leq \|\uu^\varepsilon\|_{2,\Omega} +
C\varepsilon \|\chi^\varepsilon \nabla_x (\uu^\varepsilon
-\w^\varepsilon)\|_{2,\Omega} $$
$$\leq
\|\uu^\varepsilon\|_{2,\Omega}+C\varepsilon \|\nabla_x
\uu^\varepsilon\|_{2,\Omega}+C(\varepsilon \alpha _{\mu
}^{-\frac{1}{2}})\|\chi^\varepsilon \sqrt{\alpha_\mu} \nabla_x
\w^\varepsilon\|_{2,\Omega}\leq $$
$$\|\uu^\varepsilon\|_{2,\Omega}+\frac{C\varepsilon}{\sqrt{\alpha_\lambda}}
\|\sqrt{\alpha_\lambda}\nabla_x
\w^\varepsilon\|_{2,\Omega_{s}^{\varepsilon}}+C(\varepsilon \alpha
_{\mu }^{-\frac{1}{2}})\|\sqrt{\alpha_\mu} \nabla_x
\w^\varepsilon\|_{2,\Omega_{f}^{\varepsilon}}\leq C_{0}.$$

The norm  $\|\theta^\varepsilon\|_{2,\Omega}$  we estimate with the
help of the usual Poincar\'{e}'s inequality for the difference
$(\theta^\varepsilon-\theta_{0}^\varepsilon)$ and estimate
\eqref{3.1}.

The rest of the proof is the same as for the case $\tau_0>0$.

\addtocounter{section}{1}
\setcounter{theorem}{0} \setcounter{lemma}{0}
\setcounter{proposition}{0} \setcounter{corollary}{0}
\setcounter{definition}{0} \setcounter{assumption}{0}

\begin{center} \textbf{\S4. Proof of Theorem \ref{theorem2}}
\end{center}

\textbf{4.1. Weak and two-scale limits of sequences of displacement,
temperatures and pressures.} On the strength of Theorem
\ref{theorem1}, the sequences $\{\theta^\varepsilon\}$,
$\{\nabla\theta^\varepsilon\}$,
 $\{p^\varepsilon\}$, $\{q^\varepsilon\}$, and  $\{\w^\varepsilon
\}$   are uniformly in $\varepsilon$ bounded in $L^2(\Omega_{T})$.
Hence there exist a subsequence of small parameters
$\{\varepsilon>0\}$ and functions  $\theta $, $p$, $q$, $\pi$ and
$\w$  such that
\begin{equation*}
\theta^\varepsilon \rightarrow \theta,\quad\nabla\theta^\varepsilon
\rightarrow \nabla\theta,\quad p^\varepsilon \rightarrow p,\quad
q^\varepsilon \rightarrow q, \quad
  \w^\varepsilon \rightarrow
\w
\end{equation*}
weakly in  $L^2(\Omega_T)$ as $\varepsilon\searrow 0$.

Due to Lemma \ref{Lemma2.1} there is a function $\uu^\varepsilon \in
L^\infty ((0,T);W^1_2(\Omega))$ such that $\uu^\varepsilon
=\w^\varepsilon $ in $\Omega_{s}\times (0,T)$, the family
$\{\uu^\varepsilon \}$ is uniformly in $\varepsilon$ bounded in
$L^\infty ((0,T);W^1_2(\Omega))$  and
\begin{equation*}
\uu^\varepsilon \rightarrow 0 \quad \mbox{strongly in } L^2
((0,T);W^1_2(\Omega))
\end{equation*}
as $\varepsilon \searrow 0$.

 Moreover,
\begin{equation} \label{4.1}
\chi^\varepsilon \alpha_\mu \D(\x,\w^\varepsilon) \rightarrow 0,
\end{equation}
strongly in $L^2(\Omega_T)$ as $\varepsilon \searrow 0$.

Relabeling if necessary, we assume that the sequences converge
themselves.

On the strength of Nguetseng's theorem, there exist 1-periodic in
$\y$ functions   $\Theta (\x,t,\y)$,  $P(\x,t,\y)$, $Q(\x,t,\y)$,
and $\W(\x,t,\y)$,  such that the sequences
$\{\nabla\theta^\varepsilon\}$, $\{p^\varepsilon\}$,
 $\{q^\varepsilon\}$, and  $\{\w^\varepsilon \}$ two-scale converge
 to  $\nabla_{x}\theta  +\nabla_{y}\Theta(\x,t,\y)$, $P(\x,t,\y)$,
 $Q(\x,t,\y)$,  and $\W(\x,t,\y)$ respectively.

Note that  the sequence  $\{\div_x \w^\varepsilon \}$ weakly
converges to $\div_x \w$.

The same arguments we apply for the functions  $\w_{0}^\varepsilon$
 and  $\theta_{0}^\varepsilon$:
\begin{equation*}
  \w_{0}^\varepsilon \rightarrow \w_{0}
\end{equation*}
weakly in  $L^2(\Omega_T)$  as  $\varepsilon\searrow 0$,
\begin{equation*}
\theta_{0}^\varepsilon \rightarrow \theta_{0}
\end{equation*}
weakly in   $L^2 ((0,T);W^1_2(\Omega))$  as  $\varepsilon\searrow 0$
and the sequence  $\{\div_x \w_{0}^\varepsilon \}$ weakly converges
to $\div_x \w_{0}$.\\

 \textbf{4.2. Micro- and macroscopic equations.}
\begin{lemma} \label{lemma4.1}
For all $ \x \in \Omega$ and $\y\in Y$ weak and two-scale limits of
the sequences $\{\theta^\varepsilon\}$, $\{p^\varepsilon\}$,
$\{\pi^\varepsilon\}$, $\{q^\varepsilon\}$, $\{\w^\varepsilon\}$,
$\{\nabla_x \vartheta^\varepsilon \}$  and $\{\nabla_x
\uu^\varepsilon \}$  satisfy the relations
\begin{eqnarray} \label{4.2}
&P=(1/m)\chi p,\quad Q=(1/m)\chi q=(1/m)\chi(p +\nu_0 p_{*}^{-1}
\partial p /\partial t)+\beta_{0f}\chi \theta;\\
\label{4.3} &(1-\chi)\W=0;\\
\label{4.4} & q=p +\nu_0 p_{*}^{-1} \partial p/\partial t
+m\beta_{0f}\theta;\\
\label{4.5} &(1/p_{*})p+\div_x \w =0;\\
\label{4.5.0} &((\w-\w_{0})\cdot \n(\x)=0, \quad \x \in S;\\
\label{4.5.1} & \div_y \W=0.
\end{eqnarray}
\end{lemma}

\begin{proof}

 The weak and  two-scale limiting passage in Eq.\eqref{1.4}
 yield that Eq.\eqref{4.4} and
\begin{equation} \label{4.6}
Q(\x,t,\y)=P(\x,t,\y)+\frac{\nu_0}{p_{*}}\partial
P(\x,t,\y)/\partial t+\Upsilon(\y)\theta(\x,t), \quad \y\in Y_{f}.
\end{equation}
 In order to prove  equation  \eqref{4.2}, into
Eq.\eqref{1.6} insert a test function ${\mathbf \psi}^\varepsilon
=\varepsilon {\mathbf \psi}\left(\x,t,\x / \varepsilon\right)$,
where ${\mathbf \psi}(\x,t,\y)$ is an arbitrary 1-periodic and
finite on $Y_f$ function in $\y$. Passing to the limit as
$\varepsilon \searrow 0$, we get
\begin{equation} \label{4.7}
\nabla_y Q=0, \quad \y\in Y_{f}.
\end{equation}
Next, fulfilling the two-scale limiting passage in the equality
$$(1-\chi^{\varepsilon})p^{\varepsilon}=0$$
we arrive at
$$(1-\chi )P=0$$
which along with Eqs.\eqref{4.6}--\eqref{4.7} justifies equation
\eqref{4.2}.

Eq.\eqref{4.5} and boundary condition \eqref{4.5.0}  are derived
quite similarly if we represent Eq.\eqref{1.5} in the form
\begin{equation}\label{4.8}
\frac{1}{\alpha_p}p^\varepsilon + \div_x
(\w^\varepsilon-\w_{0}^\varepsilon ) =(1-\chi^\varepsilon ) \div_x
\uu^\varepsilon -\div_x\w_{0}^\varepsilon ,
\end{equation}
multiply by an arbitrary function ${\mathbf \psi}^\varepsilon
={\mathbf \psi}(\x,t)$, integrate, and then pass to the limit as
$\varepsilon\searrow 0$. Using now in \eqref{4.8} test functions in
the form ${\mathbf \psi}^\varepsilon =\varepsilon{\mathbf
\psi}(\x,t,\x / \varepsilon )$ we obtain \eqref{4.5.1}.

  In order to prove Eq.\eqref{4.3} it is
sufficient to consider the two-scale limiting relation in
\begin{equation*}
(1-\chi ^{\varepsilon})(\w^{\varepsilon}-\uu^{\varepsilon})=0.
\end{equation*}
\end{proof}

\begin{lemma} \label{lemma4.2} For all $(\x,t) \in \Omega_{T}$
the microscopic equation
\begin{equation} \label{4.9}
\div_y \{K(\y)(\nabla_{x}\theta +\nabla_{y}\Theta) \}=0, \quad \y\in
Y
\end{equation}
holds true.  Here $K=\kappa_{0f}\chi +\kappa_{0s}(1-\chi)$.
\end{lemma}

\begin{proof}
First of all, using continuity equation \eqref{1.5} we rewrite the
heat equation in the form
\begin{equation} \label{4.10}
\alpha_\tau c^{\varepsilon}_p \frac{\partial
\theta^{\varepsilon}}{\partial t} = \div_x (\alpha^{\varepsilon}
_{\kappa} \nabla_x \theta^{\varepsilon}) -\alpha_{\theta s}
(1-\chi^{\varepsilon})\frac{\partial}{\partial t}(\div_x
\uu^{\varepsilon})+\frac{\alpha_{\theta s}}{\alpha_p}\frac{\partial
p^{\varepsilon}}{\partial t}.
\end{equation}
  Substituting now a test
function of the form ${\mathbf \psi}^\varepsilon =\varepsilon
{\mathbf \psi}\left(\x,t,\x / \varepsilon \right)$, where ${\mathbf
\psi}(\x,t,\y)$ is an arbitrary 1-periodic in $\y$ function
vanishing on the boundary $\partial \Omega$, into corresponding
integral identity, and passing to the limit as $\varepsilon \searrow
0$, we arrive at the desired microscopic relation on the cell $Y$.
\end{proof}
In the same way, using a test function independent of the fast
variable $\y/\varepsilon$, we get from Eq.\eqref{4.10}
\begin{lemma} \label{lemma4.3}
For all $(\x,t) \in \Omega_{T}$ the macroscopic equation
\begin{eqnarray}\label{4.11}
\tau_{0}\hat{c}_{p}\frac{\partial \theta}{\partial t}
-\frac{\beta_{0f}}{p_{*}}\frac{\partial p}{\partial t}= \div_x
\{(\hat{\kappa} _{0}\nabla_{x}\vartheta + \langle
K\nabla\Theta\rangle _{Y}\}
\end{eqnarray}
hold true. Here $\hat{c}_{p}=mc_{pf}+(1-m)c_{ps}, \, \hat{\kappa}
_{0}=\langle K\rangle _{Y}$.
\end{lemma}

Now we pass to the microscopic equations for the velocities in  the
liquid.
\begin{lemma} \label{lemma4.4}
Let  $\V=\chi\partial \W / \partial t$. Then
\begin{equation}\label{4.12}
\tau_{0}\rho_{f}\frac{\partial \V}{\partial t}= \mu_{1}\triangle_y
\V -\nabla_y R -\nabla_x q,\quad \y \in Y_{f}, \quad t>0;
\end{equation}
\begin{equation}\label{4.13}
\V=0,\quad \y \in \gamma ; \quad \V(\y,0)=0, \quad \y \in Y_{f}
\end{equation}
in the case $\mu_{1}>0$, and
\begin{equation}\label{4.14}
\tau_{0}\rho_{f}\frac{\partial \V}{\partial t}= -\nabla_y R
    -\nabla _{x} q , \quad \y \in Y_{f},\quad t>0;
\end{equation}
\begin{equation}\label{4.15}
 \V\cdot{\mathbf n}=0, \quad \y \in \gamma ; \quad
\V(\y,0)=0, \quad \y \in Y_{f}
\end{equation}
in the case $\mu_{1}=0$.

In Eq.\eqref{4.15} ${\mathbf n}$ is the unit normal to $\gamma$.
\end{lemma}
\begin{proof}
 Differential equations \eqref{4.12} and \eqref{4.14} follow
 as $\varepsilon\searrow 0$
 from integral equality \eqref{1.6} with the test function ${\mathbf
\psi}={\mathbf \varphi}(x\varepsilon^{-1})\cdot h({\mathbf x},t)$,
where ${\mathbf \varphi}$ is solenoidal and finite in $Y_{f}$.

 Boundary condition in  \eqref{4.13}
is the consequence of the two-scale convergence of
$\{\alpha_{\mu}^{\frac{1}{2}}\nabla_x \w^{\varepsilon}\}$ to the
function $\mu_{1}^{\frac{1}{2}}\nabla_y\W(\x,t,\y)$. On the strength
of this convergence, the function $\nabla_y \W(\x,t,\y)$ is
$L^2$-integrable in $Y$.  The boundary condition \eqref{4.15}
follows from Eq.\eqref{4.5.1}.
\end{proof}

\textbf{4.3. Homogenized equations.}

\begin{lemma} \label{lemma4.5}
Weak and strong limits   $q$ and  $\theta $  satisfy the
initial-boundary value problem
\begin{equation}\label{4.16}
\tau_{0}\hat{c_{p}}\frac{\partial \theta}{\partial t}
-\frac{\beta_{0f}}{p_{*}}\frac{\partial p}{\partial t}
 =\div_x ( B^{\theta}\cdot \nabla \vartheta ),\quad \x \in \Omega, \, t>0;
\end{equation}
\begin{equation}\label{4.17}
 \tau _{0}\theta(\x,0)=0,\quad \x \in \Omega;  \quad
 \theta(\x,0)=\theta _{0}(\x), \quad \x \in S, \, t>0,
\end{equation}
where a symmetric strictly positively defined matrix $B^{\theta}$ is
defined by formula \eqref{4.19}.
\end{lemma}
\begin{proof}
For $i=1,2,3$  we consider the periodic in $\y$ model problems
\begin{equation} \label{4.9}
\div_y \{K(\y)(\nabla_{y}\Theta_{i}+\e_{i})\}=0, \quad \y\in Y,
\end{equation}
where $({\mathbf e}_1, {\mathbf e}_2, {\mathbf e}_3)$ are the
standard Cartesian basis vectors,  and put
\begin{equation}\label{4.18}
\Theta =\sum_{i=1}^{3}(\Theta_{i} \otimes \e_{i})\cdot \nabla
_{x}\theta.
\end{equation}
Then $\Theta $ solves the  problem \eqref{4.9} and
\begin{equation}\label{4.19}
B^{\theta}=\hat{\kappa}_{0}\I+\sum_{i=1}^{3}\langle
K\nabla_{y}\Theta_{i} ^{s}\rangle _{Y}\otimes \e_{i}).
\end{equation}
All properties of the matrix $B^{\theta}$ are well known ( see
\cite{S-P}, \cite{JKO}).
\end{proof}
\begin{lemma} \label{lemma4.6}
The strong and weak limits $\theta$, $\vv=\langle \V \rangle
_{Y_{f}}$, $p$ and $q$  satisfy in $\Omega_{T}$  equation
\eqref{4.5},  Darcy's law in the form
\begin{equation}\label{4.20}
\vv=-\int_{0}^{t} B_{1}(\mu_1,t-\tau)\cdot \nabla_x q(\x,\tau )d\tau
\end{equation}
 in the case of $\tau_{0}>0$ and
$\mu_{1}>0$, Darcy's law in the form
\begin{equation}\label{4.21}
\vv=-B_{2}(\mu_1)\cdot\nabla_x q
\end{equation}
in the case of $\tau_{0}=0$  and, finally, Darcy's law in the form
\begin{equation}\label{4.22}
\vv=-\frac{1}{\tau_{0}\rho_{f}}(m\I-B_{3})\cdot\int_{0}^{t}\nabla_x
q(\x,\tau )d\tau
\end{equation}
 in the case of $\mu_{1}=0$. The
problem is supplemented by  the boundary condition
\begin{equation}\label{4.23}
(\vv-\vv_{0})\cdot \n(\x)=0, \quad \x \in S, \quad t>0,
\end{equation}
 In Eqs.
\eqref{4.20}--\eqref{4.23} $\n(\x)$ is the unit normal vector to $S$
at a point $\x\in S$, and symmetric strictly positively defined
matrices $B_{1}(\mu_1,t)$, $B_{2}(\mu_1)$, and $B_{3}$ are given in
\cite{AM}.
\end{lemma}
The proof of these statements repeats the proof of Lemma 5.8 in
\cite{AM}.

\end{document}